\documentclass[11 pt,english]{article}
\usepackage{amssymb}
\usepackage{amsmath}
\usepackage{epsfig}
\usepackage{amsthm}
\usepackage{color}

\usepackage[activeacute]{babel}
\usepackage{amsfonts} 
\usepackage{mathrsfs} 
\usepackage{amssymb, amsmath, amsfonts}
\usepackage{makeidx}
\usepackage{epsfig}
\usepackage{color}
\usepackage[T1]{fontenc}
\usepackage{graphics,graphicx,graphpap}
\usepackage[ansinew]{inputenc}
\usepackage{amsthm}
\usepackage{multirow}

\newcommand\blfootnote[1]{%
  \begingroup
  \renewcommand\thefootnote{}\footnote{#1}%
  \addtocounter{footnote}{-1}%
  \endgroup
}

\if@twoside \oddsidemargin 5pt \evensidemargin 5pt \marginparwidth 20pt
 \else \oddsidemargin 5pt \evensidemargin 5pt \marginparwidth 20pt \fi

\newcommand{\ZZ}{\mathbb{Z}}
\newcommand{\CC}{\mathbb{C}}

\newcommand{\I}{\mathcal{I}}
\newcommand{\A}{\mathcal{A}}

\newcommand{\G}{\Gamma}

\newcommand{\Circ}{\mathop{\rm Circ}}

\newcommand{\la}{\langle}
\newcommand{\ra}{\rangle}

\newtheorem{theorem}{Theorem}[section]
\newtheorem{proposition}[theorem]{Proposition}
\newtheorem{corollary}[theorem]{Corollary}
\newtheorem{lemma}[theorem]{Lemma}
\newtheorem{problem}[theorem]{Problem}

\theoremstyle{definition}

\begin{document}

\begin{center}
\Large{\textbf{On distance magic circulants of valency 6}} \\ [+4ex]
\v Stefko Miklavi\v c{\small$^{a, b, c}$}, Primo\v z \v Sparl{\small$^{a, c, d, *}$}
\\ [+2ex]
{\it \small 
$^a$University of Primorska, Institute Andrej Maru\v si\v c, Koper, Slovenia\\
$^b$University of Primorska, FAMNIT, Koper, Slovenia\\
$^c$Institute of Mathematics, Physics and Mechanics, Ljubljana, Slovenia\\ 
$^d$University of Ljubljana, Faculty of Education, Ljubljana, Slovenia\\
}
\end{center}

\blfootnote{
Email addresses: 
stefko.miklavic@upr.si (\v Stefko Miklavi\v c),
primoz.sparl@pef.uni-lj.si (Primo\v z \v Sparl)
\\
* - corresponding author
}


\hrule

\begin{abstract}
A graph $\G = (V,E)$ of order $n$ is {\em distance magic} if it admits a bijective labeling $\ell \colon V \to \{1,2, \ldots, n\}$ of its vertices for which there exists a positive integer $\kappa$ such that $\sum_{u \in N(v)} \ell(u) = \kappa$ for all vertices $v \in V$, where $N(v)$ is the neighborhood of $v$. 

A {\em circulant} is a graph admitting an automorphism cyclically permuting its vertices. In this paper we study distance magic circulants of valency $6$. We obtain some necessary and some sufficient conditions for a circulant of valency $6$ to be distance magic, thereby finding several infinite families of examples. The combined results of this paper provide a partial classification of all distance magic circulants of valency $6$. In particular, we classify distance magic circulants of valency $6$, whose order is not divisible by 12.
\end{abstract}

\hrule

\begin{quotation}
\end{quotation}

\section{Introduction}
\label{sec:Intro}

Throughout this paper all graphs are simple, undirected, finite and connected.
A graph $\G=(V,E)$ of order $n$ is {\em distance magic} if it admits a bijection $\ell \colon V \to \{1,2, \ldots, n\}$ for which there exists a positive integer $\kappa$ such that $\sum_{u \in N(v)} \ell(u) = \kappa$ for all vertices $v \in V$, where $N(v)$ is the neighborhood of $v$. In this case $\ell$ is a {\em distance magic labeling} of $\G$ and $\kappa$ is the {\em magic constant} of $\G$. Magic labelings of graphs appeared in the literature under various names. They were first defined in~\cite{Vi} as $\Sigma$-labelings and later in~\cite{MRS} as $1$-vertex magic vertex labelings. 
Following the paper by Sugeng et al. from 2009~\cite{SFMRW} the term distance magic labeling started to be used consistently. 
Since then, numerous papers on the topic appeared - see \cite{AruFroKam11} for a (now already not up to date) survey of the topic and \cite{Cic14, CicGor18, GodSinAru18, MikSpa21, MikSpa22, TiaHouHouGao21} for some recent papers on distance magic graphs. 


A simple double counting argument shows that for a $k$-regular distance magic graph $\G$ of order $n$ its magic constant is $\kappa = k(n+1)/2$. It follows that regular distance magic graphs are of even valency. Moreover, it is clear that the $4$-cycle is the only distance magic cycle. This motivated Rao~\cite{Rao08} to pose the problem of characterizing tetravalent distance magic graphs. It seems that the solution of this problem in its whole generality is beyond reach, and so Cichacz and Froncek~\cite{CicFro16} decided to study a particularly nice class of tetravalent graphs, namely tetravalent {\em circulants}, which are graphs admitting an automorphism cyclically permuting its vertices. They obtained a partial classification of these graphs (see~\cite[Theorem 16]{CicFro16}). A complete classification of distance magic tetravalent circulants was recently achieved in~\cite{MikSpa21}.

With the classification of distance magic tetravalent circulants completed, there are now two obvious directions for further research: the first one is to study other nice families of distance magic tetravalent graphs, and the second one is to study distance magic circulants of the next admissible valency, namely $6$. In this paper we initiate the latter direction of research. 

As in the case of the investigation of distance magic tetravalent circulants in~\cite{MikSpa21}, the main tool for studying distance magic circulants of valency $6$ in this paper will be the irreducible characters of cyclic groups. It turns out that the so-called admissible irreducible characters of the cyclic group $\ZZ_n$ (see Section~\ref{sec:eig0} for the definition and details), which play the key role in our investigation of $6$-valent distance magic circulants of order $n$, come in three different types. In this paper we obtain a complete classification of $6$-valent distance magic circulants for which all admissible irreducible characters are of the same type (see Theorem~\ref{the:type1class}, Theorem~\ref{the:type2class} and Proposition~\ref{pro:type3}). We further show that no $6$-valent distance magic circulant can have admissible irreducible characters of all three types and that at most two combinations of two types are possible. For such examples we obtain some partial results. We also obtain infinite families of nontrivial examples of distance magic circulants of valency $6$ (in the sense that they are not the rather ``obvious'' examples of lexicographic products of a prism or a M\"obius ladder by $2K_1$ or of a cycle by $3K_1$ - see Section~\ref{sec:Prelim}). Finally, the combined results of this paper provide a complete classification of distance magic circulants of valency 6, whose order is not divisible by $12$ (see Theorem~\ref{the:class_part}).

\section{Preliminaries}
\label{sec:Prelim}

We first set some standard notation that will be used throughout the paper. For an integer $n$ we let $\ZZ_n$ denote the ring of residue classes modulo $n$ (and at the same time the cyclic group of order $n$) and we let $\ZZ_n^*$ be its group of units. For a subset $S$ of $\ZZ_n$ and an element $q \in \ZZ_n^*$ we let $qS = \{qs \colon s \in S\}$. For $a,b \in \ZZ_n$ we let $a = \pm b$ mean that $a$ is equal to $b$ or $-b$, while we let $a \neq \pm b$ mean that $a$ is not equal to $b$ nor to $-b$. 

As stated in the Introduction we will be studying distance magic circulants of valency $6$ in this paper. Recall that a circulant is a Cayley graph of a cyclic group. More precisely, for an integer $n$ and a subset $S \subset \ZZ_n$ with $S = -S$ and $0 \notin S$ the {\em circulant} $\Circ(n ; S)$ is the graph with vertex-set $\ZZ_n$ in which distinct vertices $i$ and $j$ are adjacent if and only if $j-i \in S$ (where the computation is done modulo $n$). The comments from the previous section imply that if a circulant $\Circ(n ; S)$ with $|S| = 6$ is distance magic, then its magic constant is
\begin{equation}
\label{eq:k}
	\kappa = 3(n+1).
\end{equation}

As mentioned in the Introduction there are some rather obvious examples of distance magic circulants of valency $6$. To describe them we first recall the definition of a lexicographic product of graphs. For graphs $\G_1 = (V_1, E_1)$ and $\G_2 = (V_2, E_2)$ the {\em lexicographic product} $\G_1[\G_2]$ of $\G_1$ by $\G_2$ is the graph with vertex-set $V_1 \times V_2$ in which vertices $(u_1,u_2)$ and $(v_1,v_2)$ are adjacent whenever either $u_1 = v_1$ and $u_2 \sim v_2$ in $\G_2$, or $u_1 \sim v_1$ in $\G_1$. It is well known and easy to verify that the well-known {\em M\"obius ladder} $\mathrm{Ml}_m$ of order $2m$, where $m \geq 3$, is the cubic circulant $\Circ(2m ; \{\pm 1, m\})$. Similarly, for an odd $m \geq 3$ the {\em prism} $\mathrm{Pr}_m$ of order $2m$ is the circulant $\Circ(2m ; \{\pm 2, m\})$ (we point out that prisms are defined also for $m$ even, but in that case they are not circulants, and so in this paper we will only consider $\mathrm{Pr}_m$ for $m$ odd). It is straightforward to verify that 
$$
	\mathrm{Ml}_m[2K_1] \cong \Circ(4m ; \{\pm 1, \pm m, \pm (2m-1)\})\ \text{and}\ \mathrm{Pr}_m[2K_1] \cong \Circ(4m ; \{\pm 2, \pm m, \pm (2m-2)\}).
$$
That $\mathrm{Ml}_m[2K_1]$ for $m \geq 3$ and $\mathrm{Pr}_m[2K_1]$ for an odd $m \geq 3$ are examples of distance magic circulants of valency $6$ is easy to see, since one can simply assign the labels $i$ and $4m+1-i$, $1 \leq i \leq 2m$, to the $2m$ pairs corresponding to $2K_1$ (but see also~\cite[Lemma~2]{MRS}). It is also straightforward to verify that for each $m \geq 3$ we have that
$$
	C_m[3K_1] \cong \Circ(3m ; \{\pm 1, \pm (m-1), \pm (m+1)\}),
$$
where $C_m$ is the cycle of length $m$. By~\cite[Theorem~3.5]{ShaAliSim09} the graph $C_m[3K_1]$ is distance magic if and only if $m$ is either odd or is divisible by $4$. In this paper we will say that the distance magic circulants $\mathrm{Ml}_m[2K_1]$, $m \geq 3$, $\mathrm{Pr}_m[2K_1]$, $m \geq 3$ odd, and $C_m[3K_1]$, $m \geq 3$ odd or divisible by $4$, are {\em trivial} examples of $6$-valent distance magic circulants. Our aim is to investigate the nontrivial ones. 

As was pointed out in~\cite{MikSpa21} the property of being distance magic is very nicely characterized by eigenvalues and eigenvectors of the adjacency matrix when it comes to regular graphs. Since this result will play an important role in our arguments, we restate it here for ease of reference.

\begin{lemma}[\cite{MikSpa21}]
\label{le:eig0}
Let $\G=(V,E)$ be a regular graph of order $n$ and even valency. Then $\G$ is distance magic if and only if $0$ is an eigenvalue of the adjacency matrix of $\G$ and there exists an eigenvector for the eigenvalue $0$ with the property that a certain permutation of its entries results in the arithmetic sequence 
\begin{equation}
\label{eq:seq}
	\frac{1-n}{2}, \frac{3-n}{2}, \frac{5-n}{2}, \ldots, \frac{n-3}{2}, \frac{n-1}{2}.
\end{equation}
In particular, if $\G$ is distance magic then $0$ is an eigenvalue for the adjacency matrix of $\G$ and there exists a corresponding eigenvector all of whose entries are pairwise distinct. 
\end{lemma}

It was also pointed out in~\cite{MikSpa21} that when studying the property of being distance magic for Cayley graphs of abelian groups, group characters can be very useful. The situation is particularly nice in the case of cyclic groups since its characters are very easy to describe. We give here only the necessary ingredients for our arguments, but refer the reader to~\cite{MikSpa21} for details. 

It is well known that the group $\hat{\ZZ}_n$ of irreducible characters for the cyclic group $\ZZ_n$ consists of the $n$ homomorphisms $\chi_j$, $0 \leq j < n$, from $\ZZ_n$ to the unit circle in the complex field $\CC$ where $\chi_j$ maps according to the rule
\begin{equation}
	\label{eq:characters}
  \chi_j(x) =\Big( e^{\frac{2 \pi \mathbf{i}}{n}} \Big)^{jx} = \cos{\Bigg( \frac{2 \pi  j x}{n} \Bigg)} + \mathbf{i} \sin{\Bigg( \frac{2 \pi  j x}{n} \Bigg)}.
\end{equation}
What is more, by~\cite[Lemma~9.2, p. 246]{G} the spectrum (that is, the set of all eigenvalues of the adjacency matrix) of the circulant $\Circ(n ; S)$ is given by 
\begin{equation}
\label{eq:spec}
	\left\{\chi_j(S) \colon 0 \leq j \leq n-1\right\},
\end{equation}
where $\chi_j(S) = \sum_{s \in S}\chi_j(s)$. Moreover, letting $\mathbf{w}_\chi$ denote the column vector with its entries indexed by $g \in \ZZ_n$ such that the $g$-entry of $\mathbf{w}_\chi$ is equal to $\chi(g)$, the $n$ vectors $\mathbf{w}_{\chi_j}$, $0 \leq j < n$, are a complete set of eigenvectors for the adjacency matrix of $\Circ(n ; S)$ and $\mathbf{w}_{\chi_j}$ corresponds to the eigenvalue $\chi_j(S)$.

\section{The eigenvalue $0$}
\label{sec:eig0}

By Lemma~\ref{le:eig0} and the discussion at the end of the previous section the only candidates for distance magic circulants $\G = \Circ(n ; S)$ are the ones for which $\chi_j(S) = 0$ for at least one $j$, $0 \leq j < n$. This motivates the following definition. For a positive integer $n$ we let $\I_n = \{0,1,\ldots , n-1\}$ and for a subset $S \subset \ZZ_n$ with $S = -S$ and $0 \notin S$ we let 
$$
	\A_n(S) = \{j \in \I_n \colon \chi_j(S) = 0\}.
$$
We say that $\A_n(S)$ is the set of all {\em admissible} elements of $\I_n$ with respect to $S$ and also say that the corresponding irreducible characters $\chi_j$ are {\em admissible} for $S$. When the set $S$ is clear from the context we simply write $\A_n$ and say that $j$ (or $\chi_j$) is admissible. The following straightforward result will be useful.

\begin{lemma}
\label{le:common_divisor}
Let $\G = \Circ(n; S)$ be a circulant. If there exists a divisor $d$ of $n$ with $1 < d < n$ dividing each $j \in \A_n(S)$, then $\G$ is not distance magic.
\end{lemma}

\begin{proof}
Suppose such a divisor $d$ exists and take any $j \in \A_n(S)$. By~\eqref{eq:characters} we have that $\chi_j(0) = 1 = \chi_j(n/d)$. Recall that the eigenvectors $\mathbf{w}_{\chi_j}$, where $j$ runs through the whole set $\A_n(S)$, constitute a basis for the eigenspace corresponding to the eigenvalue $0$ of the adjacency matrix of $\G$. Therefore, each eigenvector for the eigenvalue $0$ of $\G$ has equal entries at $0$ and $n/d$, and so $\G$ cannot be distance magic by Lemma~\ref{le:eig0}.
\end{proof}

Now, let $n \geq 7$ be an integer, let $S = \{\pm a, \pm b, \pm c\} \subset \ZZ_n$ be such that $|S| = 6$ and let $\G = \Circ(n; S)$. Suppose that the set $\A_n(S)$ is nonempty (by Lemma~\ref{le:eig0} this must be the case if $\G$ is distance magic). By the discussion of the previous section and \eqref{eq:characters} this holds if and only if 
\begin{equation}
\label{eq:cos_sum}
	\cos{\Bigg( \frac{2 \pi  j a}{n} \Bigg)} + \cos{\Bigg( \frac{2 \pi  j b}{n} \Bigg)} + \cos{\Bigg( \frac{2 \pi  j c}{n} \Bigg)} = 0.
\end{equation}
Observe that each of $\frac{2ja}{n}$, $\frac{2jb}{n}$ and $\frac{2jc}{n}$ is a rational number. Moreover, since $\cos(x) = \cos(-x)$ for all $x$ we can (at least for the purposes of considering the possible solutions of~\eqref{eq:cos_sum}) assume that $0 \leq \frac{2ja}{n}, \frac{2jb}{n}, \frac{2jc}{n} \leq 1$. The problem, posed in 1944 by H.~S.~M.~Coxeter, of determining all rational solutions of the equation
\begin{equation}
\label{eq:Cox}
	\cos(r_1\pi) + \cos(r_2\pi) + \cos(r_3\pi) = 0,\quad 0 \leq r_1 \leq r_2 \leq r_3 \leq 1,
\end{equation}
was solved in 1946 by W.~J.~R.~Crosby~\cite{CoxCro46}. It was proved that, except for a pair of ``symmetric'' exceptions, the only solutions of~\eqref{eq:Cox} are those that belong to two infinite families of ``obvious'' triples $(r_1, r_2, r_3)$, namely 
\begin{equation}
\label{eq:sol1}
	0 \leq r_1 \leq \frac{1}{2},\quad r_2 = \frac{1}{2},\quad r_3 = 1 - r_1,
\end{equation}
and
\begin{equation}
\label{eq:sol2}
	0 \leq r_1 \leq \frac{1}{3},\quad r_2 = \frac{2}{3}-r_1,\quad r_3 = \frac{2}{3} + r_1.
\end{equation}
The only two exceptions are 
\begin{equation}
\label{eq:sol3}
	r_1 = \frac{1}{5},\ r_2 = \frac{3}{5},\ r_3 = \frac{2}{3}\quad \text{and}\quad r_1 = \frac{1}{3},\ r_2 = \frac{2}{5},\ r_3 = \frac{4}{5}.
\end{equation}

It is clear that no triple $(r_1,r_2,r_3)$ of rational numbers with $0 \leq r_1 \leq r_2 \leq r_3 \leq 1$ which satisfies any of the two possibilities from~\eqref{eq:sol3} satisfies~\eqref{eq:sol1} or~\eqref{eq:sol2}. Moreover, the only triple $(r_1,r_2,r_3)$ which satisfies both~\eqref{eq:sol1} and~\eqref{eq:sol2} is $(\frac{1}{6}, \frac{1}{2}, \frac{5}{6})$. We can thus introduce the following terminology. For a given integer $n \geq 7$ and a subset $S = \{\pm a, \pm b, \pm c\} \subset \ZZ_n$ of size $6$ we say that a $j \in \A_n(S)$ (as well as the corresponding character $\chi_j$) is of {\em type 1}, {\em type 2} or {\em type 3}, respectively, if the corresponding solution of~\eqref{eq:Cox} is of type~\eqref{eq:sol1}, \eqref{eq:sol2} or~\eqref{eq:sol3}, respectively. Except for the one above mentioned exception (where a $j \in \A_n(S)$ can be of types 1 and 2) each $j \in \A_n(S)$ thus has a unique type. We point out that whenever we will say that {\em all $j \in \A_n(S)$ are of type~1} (type~2, respectively) we will mean that each $j \in \A_n(S)$ is of type~1 (type~2, respectively) but we do allow the above mentioned possibility that some $j \in \A_n(S)$ are at the same time of type~2 (type~1, respectively). In the following sections we analyze each of the three types separately.

\section{Type~1}
\label{sec:type1}

In this section we consider admissible $j \in \A_n(S)$ of type 1 and classify the distance magic circulants of valency $6$ for which all admissible characters are of type~1. Before stating our first result we introduce the following terminology that will be used throughout the rest of the paper. For an integer $m$ and a prime $p$ we let the {\em $p$-part} of $m$ be $p^t$, where $t$ is the largest integer such that $m$ is divisible by $p^t$.

\begin{lemma}
\label{le:type1}
Let $n \geq 7$ be an integer and let $S = \{\pm a, \pm b, \pm c\} \subset \ZZ_n$ be such that $|S| = 6$ and $\la S \ra = \ZZ_n$. Suppose at least one $j \in \A_n(S)$ of type~1 exists. Then $n = 4n_0$ for an integer $n_0 \geq 2$. Moreover, for each such $j$ there are $s_1, s_2, s_3 \in S$ with $\{\pm s_1, \pm s_2, \pm s_3 \} = S$ such that $s_1 + s_3$ is even and 
\begin{equation}
\label{eq:type1}
	js_2 = n_0(1+2k_1) \quad \text{and}\quad j(s_1 + s_3) = 2n_0(1 + 2k_2)
\end{equation}
for some integers $k_1, k_2$, where these two equations are to be read within the ring of integers $\ZZ$.  
\end{lemma}

\begin{proof}
Let $j \in \A_n(S)$ be such that \eqref{eq:cos_sum} holds and suppose that the corresponding solution of~\eqref{eq:Cox} is of type~\eqref{eq:sol1}. There thus is some $s_2 \in \{a,b,c\}$ such that $\cos(2\pi js_2/n) = 0$ (that is, $2js_2/n$ corresponds to $r_2$ from~\eqref{eq:sol1}). Moreover, letting $s_1, s_3 \in S \setminus \{\pm s_2\}$ be such that $2js_1/n$ and $2js_3/n$ correspond to $r_1$ and $r_3$ we therefore have that
$$
	\frac{2\pi js_2}{n} = \frac{\pi}{2} + k_1\pi \quad \text{and}\quad \frac{2\pi js_1}{n} = \pi - \frac{2\pi js_3}{n} + 2k_2\pi	
$$
for some integers $k_1$ and $k_2$. Rearranging we obtain
$$
	4js_2 = n(1+2k_1)\quad \text{and}\quad 2j(s_1+s_3) = n(1 + 2k_2),
$$
and so the first of these two equations implies $n = 4n_0$ for some $n_0 \geq 2$. The above two equations then transform into~\eqref{eq:type1}. Note that~\eqref{eq:type1} implies that the 2-part of $j(s_1+s_2)$ is larger than the 2-part of $js_2$, and so $s_1+s_3$ must be even. 
\end{proof}

We can now characterize the only candidates for distance magic $6$-valent circulants for which all admissible characters are of type~1. 

\begin{proposition}
\label{pro:type1}
Let $n \geq 7$ be an integer and let $S = \{\pm a, \pm b, \pm c\} \subset \ZZ_n$ be such that $|S| = 6$ and $\la S \ra = \ZZ_n$. If $\G = \Circ(n; S)$ is distance magic and all $j \in \A_n(S)$ are of type~1, then $n = 4n_0$ for some $n_0 \geq 2$. Moreover, either $\G$ is trivial in the sense that there exists some $q \in \ZZ_n^*$ such that $qS = \{\pm 1, \pm n_0, \pm (2n_0-1)\}$ or $qS = \{\pm 2, \pm n_0, \pm(2n_0-2)\}$ in which case $\G \cong \mathrm{Ml}_{n_0}[2K_1]$ or $\G \cong \mathrm{Pr}_{n_0}[2K_1]$, or one of the following holds:
\begin{itemize}
\item $n_0 = dd'$ for odd and coprime integers $d, d'$ with $1 < d < d' < n_0$ and there exists $q \in \ZZ_n^*$ such that $qS = \{\pm 2, \pm n_0, \pm c'\}$, where $1 < c' < n$ is the unique solution of the system of congruences
\begin{equation}
\label{eq:T1sys1}
	\begin{array}{r@{\,}c@{\,}r@{}l}
	c' & \equiv & 0 & \pmod{4}\\
	c' & \equiv & 2 & \pmod{d}\\
	c' & \equiv & -2 & \pmod{d'}.
	\end{array}
\end{equation}
\item $n_0 = dd'd''$ for odd and coprime integers $d, d', d''$ with $1 \leq d < d' < d'' < n_0$ and there exists $q \in \ZZ_n^*$ such that $qS = \{\pm d, \pm b', \pm c'\}$, where $1 < b', c' < n$ are the unique solutions of the systems of congruences
\begin{equation}
\label{eq:T1sys2}
\begin{array}{r@{\,}c@{\,}c@{}lcr@{\,}c@{\,}c@{}l}
	b' & \equiv & 2-d & \pmod{4} &  &  c' & \equiv & 2-d & \pmod{4}\\
	b' & \equiv & 0 & \pmod{d'} & \text{and} & c' & \equiv & -b' & \pmod{d}\\
	b' & \equiv & -d & \pmod{d''}&   & c' & \equiv & -d & \pmod{d'}\\
	   &          &    &               &   & c' & \equiv & 0 & \pmod{d''}
\end{array}
\end{equation}
subject to the condition that $1 < b' < n$ is the smallest solution of the left system with $\gcd(b',d) = 1$.
\end{itemize}
\end{proposition}

\begin{proof}
Suppose that $\G$ is distance magic and all $j \in \A_n(S)$ (abbreviated to $\A_n$ in the rest of this proof) are of type~1. Lemma~\ref{le:type1} thus imples that $n = 4n_0$ for some $n_0 \geq 2$, while Lemma~\ref{le:common_divisor} implies that there is no divisor $d > 1$ of $n$ dividing all $j \in \A_n$. Suppose that $\G$ is none of the two trivial examples of lexicographic products from the statement of the proposition. We proceed by proving a series of claims. Before stating and proving the first of them we make the agreement that 
throughout the rest of the proof for a $j \in \A_n$ we let $s_1, s_2, s_3 \in S$ and $k_1, k_2 \in \ZZ$ be as in~\eqref{eq:type1} and similarly we let $s'_1, s'_2, s'_3 \in S$ and $k'_1, k'_2 \in \ZZ$ be as in~\eqref{eq:type1} for a $j' \in \A_n$.\smallskip


\noindent
{\sc Claim 1:} For any prime divisor of $n_0$ precisely one of $a$, $b$ and $c$ is divisible by the whole $p$-part of $n_0$, while the remaining two are coprime to $p$.\\
Let $p$ be a prime divisor of $n_0$ and let $p^t$ be the $p$-part of $n_0$. By Lemma~\ref{le:common_divisor} there is a $j \in \A_n$, which is coprime to $p$. Then Lemma~\ref{le:type1} implies that $p^t$ divides $s_2$ and $s_1 + s_3$. If $p$ divided any of $s_1$ and $s_3$, it would thus have to divide both, which is impossible as $S$ generates $\ZZ_n$. 
\smallskip

\noindent
{\sc Claim 2:} All admissible $j \in \A_n$ are odd.\\
Suppose to the contrary that this is not the case. By Lemma~\ref{le:common_divisor} we then have $j, j' \in \A_n$ such that $j$ is even, while $j'$ is odd. Denote the $2$-part of $n_0$ by $2^t$. Since $j$ is even,~\eqref{eq:type1} implies that $t \geq 1$ and the 2-part of $s_2$ is smaller than $2^t$. On the other hand, as $j'$ is odd, the 2-part of $s'_2$ equals $2^t$ (which in turn implies that $s'_2$ is even), and so $s'_2$ and $n-s'_2$ are both divisible by $2^t$, implying that $s_2 \neq \pm s'_2$. It follows that $s'_2 \in \{\pm s_1, \pm s_3\}$. Since $s_1 + s_3$ is even by Lemma~\ref{le:type1}, we thus find that $s_1$ and $s_3$ are both even, contradicting Claim~1.
\smallskip

\noindent
{\sc Claim 3:} $n_0$ is odd.\\
By way of contradiction suppose $n_0$ is even. Since all $j \in \A_n$ are odd, Lemma~\ref{le:type1} implies that for each $j \in \A_n$ the corresponding $s_2$ is even. Claim~1 thus implies that one of $a$, $b$ and $c$ is even, while the other two are odd. Consequently, for all $j \in \A_n$ the corresponding $s_2$ is the same (modulo multiplication by $-1$). Denote it by $s$. Then~\eqref{eq:characters} and~\eqref{eq:type1} imply that $\chi_j(4s) = 1$ for all $j \in \A_n$, and so Lemma~\ref{le:eig0} forces $4s = 0$ in $\ZZ_n$. Since $s \neq -s$ it follows that $s = \pm n_0$. Since for each $j \in \A_n$ the corresponding $s_1$ and $s_3$ are both odd and $2n_0$ is divisible by $4$, $s_1$ and $s_3$ have different odd remainders modulo $4$. Moreover, for any $j, j' \in \A_n$ we must have that $s'_1 + s'_3 = \pm (s_1 + s_3)$, and so letting $h$ be this sum (modulo multiplication by $-1$) we have that $\chi_j(2h) = 1$ for all $j \in \A_n$. Thus $h = 2n_0$ (since $s_1 + s_3 \neq 0$). Multiplying by a suitable $q \in \ZZ_n^*$ we get $qS = \{\pm 1, \pm n_0, \pm (2n_0 - 1)\}$ (note that $s_1$ and $s_3$ are coprime to $n$ by Claim~1), contradicting our hypothesis that $\G$ is not trivial. 
\smallskip

\noindent
We now consider two different cases, depending on whether two of the elements from $\{a,b,c\}$ are coprime to $n_0$ or not. Note that by Claim~1 at least one of them is not coprime to $n_0$.
\smallskip

\noindent
{\sc Case 1:} Two elements of $\{a,b,c\}$ are coprime to $n_0$: \\
By Claim~1 one of $a$, $b$ and $c$, with no loss of generality assume it is $b$, is divisible by $n_0$. Since $b \neq -b$, it thus follows that $b = \pm n_0$. Let $p$ be any prime divisor of $n_0$. By Lemma~\ref{le:common_divisor} there is some $j \in \A_n$ which is coprime to $p$. Lemma~\ref{le:type1} then implies that the corresponding $s_2$ is divisible by $p$, and so $s_2 = \pm b$. Moreover, the $p$-part of $n_0$ divides $s_1 + s_3$. If for each pair of primes $p$ and $p'$ dividing $n_0$ there are corresponding $j,j' \in \A_n$ with $p \nmid j$ and $p'\nmid j'$ such that $s'_1 + s'_3 = \pm (s_1 + s_3)$, then $n_0$ must divide this common sum (say $\pm(s_1 + s_3)$), and so since by~Lemma~\ref{le:type1} this sum is even, it is $2n_0$. But then $\G$ is trivial, a contradiction. It thus follows that there are two distinct primes $p,p'$ dividing $n_0$ with corresponding $j,j' \in \A_n$ such that $s'_1+s'_3 = \pm(s_1 - s_3)$. Hence, $s_1$ and $s_3$ must be even (otherwise one of $s_1 \pm s_3$ is divisible by $4$, contradicting~\eqref{eq:type1} and Claim~3) and in addition one of them is divisible by $4$ while the other is not. Multiplying by an appropriate $q \in \ZZ_n^*$ if necessary we can thus assume that $a = 2$. By what we have just shown $c$ is divisible by $4$ and for each prime divisor $p$ of $n_0$ the $p$-part of $n_0$ divides precisely one of $c-2$ and $c+2$. Let $d$ be the product of all the $p$-parts of $n_0$ dividing $c-2$ and let $d'$ be the product of all the $p$-parts of $n_0$ dividing $c+2$. Of course, $n_0 = dd'$. Replacing the roles of $c$ and $-c$ if necessary we can assume that $d < d'$. By the Chinese remainder theorem the fact that $4$, $d$ and $d'$ are pairwise coprime implies that the conditions $4 \mid c$, $d \mid (c-2)$ and $d' \mid (c+2)$ determine $c$ uniquely (recall that $n = 4dd'$).
\smallskip

\noindent
{\sc Case 2:} At most one element of $\{a,b,c\}$ is coprime to $n_0$:\\
Let $d = \gcd(a,n_0)$, $d' = \gcd(b,n_0)$ and $d'' = \gcd(c,n_0)$, where with no loss of generality we assume $d \leq d' \leq d''$. Claim~1 implies that $d$, $d'$ and $d''$ are pairwise coprime and $n_0 = dd'd''$. This also implies that $d < d' < d''$ and in particular $d' > 1$ (but we may have $d = 1$). Let $p$ be any prime dividing $d'$. By Lemma~\ref{le:common_divisor} there exists a $j \in \A_n$ coprime to $p$, and so Lemma~\ref{le:type1} implies that the corresponding $s_2$ is $\pm b$. Then Claim~3 and \eqref{eq:type1} imply that $b$ is odd and that the whole $p$-part of $n_0$ divides $c + \delta a$, where $\delta \in \{-1,1\}$ is the unique element such that $c + \delta a \equiv 2 \pmod{4}$. Consequently, the whole $d'$ divides $c + \delta a$. In a completely analogous way we see that $c$ is also odd and that for the unique $\delta' \in \{-1, 1\}$ such that $b + \delta' a \equiv 2 \pmod{4}$ the whole $d''$ divides $b + \delta' a$. As $b$ is odd this implies that $a$ must also be odd. Multiplying by an appropriate $q$ if necessary we can thus assume that $a = d$. Moreover, by definition of $d'$ and $d''$ we have that $b = \xi' d'$ and $c = \xi'' d''$ for appropriate $\xi'$ coprime to $2dd''$ and $\xi''$ coprime to $2dd'$, where $1 \leq \xi' < 4dd''$ and $1 \leq \xi'' < 4dd'$.
Replacing the roles of $b$ and $-b$ and of $c$ and $-c$ if necessary we can assume that $\delta = \delta' = 1$. Therefore, $b$ is a solution of the system of congruences from the left-hand side of~\eqref{eq:T1sys2}. If $d = 1$ then $b$ is uniquely determined (since $n = 4d'd''$) and is coprime to $d$. Moreover, the condition that $c \equiv -b \pmod{d}$ is trivial, and so $c$ is the (unique) solution of the system of congruences from the right-hand side of~\eqref{eq:T1sys2}. Suppose finally that $d > 1$. Just as above we find that the whole $d$ must divide $c+b$ (note that since $b+a$ and $c+a$ are both even but not divisible by $4$, $c-b$ is divisible by $4$), and so $c \equiv -b \pmod{d}$. Therefore, $c$ satisfies the system of congruences from the right-hand side of~\ref{eq:T1sys2}. Note that once $b$ is fixed $c$ is completely determined by~\eqref{eq:T1sys2}. To see that we can assume $b$ is the smallest solution of the system from the left-hand side of~\eqref{eq:T1sys2} such that $\gcd(b,d) = 1$, suppose $b_1$ and $b_2$, $1 < b_1, b_2 < n$ are two solutions of that system with $\gcd(b_1,d) = \gcd(b_2,d) = 1$. Then $b_2 - b_1$ is divisible by $4d'd''$, say $b_2 = b_1 + 4\xi d'd''$ with $0 < \xi < d$. Let $c_1$ and $c_2$, $1 < c_1,c_2 < n$, be the corresponding unique solutions of the system on the right-hand side of~\eqref{eq:T1sys2}. Since $\gcd(b_1,d) = 1$, there exists some $\eta$, $1 \leq \eta < d$, such that $\eta b_1 \equiv 1 \pmod{d}$. Set $q = 4\eta\xi d'd''+1$ and note that $q$ is coprime to $4d'd''$ and that $b_2 = qb_1$ in $\ZZ_n$. Since $\gcd(b_2,d) = 1$ this in fact implies that $q \in \ZZ_n^*$. Clearly, $qd \equiv d \pmod{n}$. It is now easy to see that $c_2 = qc_1$ (as $qc_1$ satisfies all four conditions for $c_2$ in~\eqref{eq:T1sys2} with $b = qb_1$). Therefore, all of the possible solutions for $b$ and $c$ are equivalent up to multiplication by a suitable $q \in \ZZ_n^*$ preserving $d$.
\end{proof}

As an example let us determine (by using Proposition~\ref{pro:type1}) the only candidates for possible distance magic circulants of the form $\Circ(1540; \{\pm a, \pm b, \pm c\})$, where $a < b < c < 770$. We of course have the lexicographic products $\Circ(1540; \{\pm 1, \pm 385, \pm 769\})$ and $\Circ(1540; \{\pm 2, \pm 385, \pm 768\})$. Since $1540 = 4\cdot 5 \cdot 7 \cdot 11$, we have all of the possibilities from the above proof. Up to multiplication by an appropriate $q \in \ZZ_{1540}^*$ we have the following possibilities. For the possibility when two of $\{a,b,c\}$ are coprime to $n_0$ we get $\Circ(1540; \{\pm 2, \pm 152, \pm 385\})$ (when $d = 5$), $\Circ(1540; \{\pm 2, \pm 385, \pm 548\})$ (when $d = 7$) and $\Circ(1540; \{\pm 2, \pm 68, \pm 385\})$ (when $d = 11$). For the other possibility we get three with $d = 1$, namely $\Circ(1540; \{\pm 1, \pm 155, \pm 231\})$ (when $d' = 5$), $\Circ(1540; \{\pm 1, \pm 329, \pm 715\})$ (when $d' = 7$) and $\Circ(1540; \{\pm 1, \pm 209, \pm 595\})$ (when $d' = 11$), while there is a unique one for $d = 5$ (and thus $d' = 7$ and $d'' = 11$), namely $\Circ(1540; \{\pm 5, \pm 413, \pm 737\})$. As we show in the next two lemmas, all of these graphs are in fact distance magic. 
\medskip

We already mentioned in Section~\ref{sec:Prelim} that the trivial examples (the two lexicographic products) from Proposition~\ref{pro:type1} are indeed distance magic. For the nontrivial examples we first show that the graphs corresponding to the first item from Proposition~\ref{pro:type1} are distance magic.

\begin{lemma}
\label{le:T1DMcase1}
Let $d,d' > 1$ be odd and coprime integers and let $n_0 = dd'$ and $n = 4n_0$. Let $c$ with $1 < c < n$ be the unique solution of the system of congruences~\eqref{eq:T1sys1}. Then the circulant $\Circ(n ; \{\pm 2, \pm n_0, \pm c\})$ is distance magic.
\end{lemma}

\begin{proof}
Let $c_0 = c/2$ and note that $c_0$ is even by~\eqref{eq:T1sys1}. Consider the subgraph $\Circ(n; \{\pm 2, \pm c\})$ of $\G = \Circ(n ; \{\pm 2, \pm n_0, \pm c\})$. Since $c$ is even, it is not connected. In fact, it consists of two isomorphic copies of the graph $\Delta = \Circ(2n_0 ; \{\pm 1, \pm c_0\})$. It follows from~\eqref{eq:T1sys1} that $c_0^2-1$ is divisible by $n_0$. Therefore, the fact that $c_0$ is even shows that $\Delta$ satisfies the conditions of~\cite[Theorem~1.1]{MikSpa21}, which thus implies that $\Delta$ is distance magic. Let $\ell_\Delta \colon \ZZ_{2n_0} \to \{1,2,\ldots , 2n_0\}$ be the corresponding distance magic labeling from the proof of~\cite[Theorem~1.1]{MikSpa21}. With a slight abuse of notation we now define a labeling $\ell \colon \ZZ_{n} \to \{1,2,\ldots , n\}$ by the following rule:
\begin{equation}
\label{eq:ell1}
	\ell (x) = \left\{\begin{array}{ccc}
		\ell_\Delta(x/2) + 2n_0 & : & i \equiv 0 \pmod{4},\\
		\ell_\Delta((x-1)/2)       & : & i \equiv 1 \pmod{4},\\
		\ell_\Delta(x/2)            & : & i \equiv 2 \pmod{4},\\
		\ell_\Delta((x-1)/2) +2n_0 & : & i \equiv 3 \pmod{4}.\end{array}\right.
\end{equation}
Using the properties of $\ell_\Delta$ we can show that $\ell$ is a distance magic labeling for $\G$. To see this we first recall that the labeling $\ell_\Delta$ assigns the labels from $\{1,2,\ldots , n_0\}$ to all the vertices $y \in \{0,1,\ldots, 2n_0-1\}$ with $y$ even. This guarantees that $\ell$ is a bijection. The other useful property of $\ell_\Delta$ is that $\ell_\Delta(y) + \ell_\Delta(y+n_0) = 2n_0+1$ for all $y \in \{0,1,\ldots , n_0-1\}$. Together with the observation that, for any $x \in \ZZ_{n}$, either precisely one of $x + n_0$ and $x - n_0$ is congruent to $1$ modulo $4$ (and the other is congruent to $3$ modulo $4$), or precisely one of them is divisible by $4$ (and the other is congruent to $2$ modulo $4$), this clearly shows that 
\begin{equation}
\label{eq:aux1}
	\ell(x+n_0) + \ell(x-n_0) = 2n_0+1 + 2n_0 = n+1.
\end{equation} 
Note that $x + 2 \equiv x - 2 \pmod{4}$, $x + c \equiv x - c \pmod{4}$ and $x + 2 \equiv x + c + 2 \pmod{4}$. Therefore, by definition of $\ell$ from~\eqref{eq:ell1} and since $\ell_\Delta$ is a distance magic labeling for $\Delta$ (whose magic constant is $2(2n_0+1)$), we thus find that 
$$
	\ell(x+2) + \ell(x-2) + \ell(x+c) + \ell(x-c) = 2(2n_0+1) + 2\cdot 2n_0 = 2(n+1)
$$
holds for all $x$. By~\eqref{eq:k} and~\eqref{eq:aux1} this proves that $\ell$ is indeed a distance magic labeling for $\G$.
\end{proof}

We next show that the graphs corresponding to the second item from Proposition~\ref{pro:type1} are also distance magic.

\begin{lemma}
\label{le:T1DMcase2}
Let $d,d',d''$, where $1 \leq d < d' < d''$ be odd and coprime integers and let $n_0 = dd'd''$ and $n = 4n_0$. Let $S = \{\pm d, \pm b, \pm c\}$, where $b$ and $c$ are the solutions of~\eqref{eq:T1sys2} with $\gcd(b,d) = 1$. Then the circulant $\Circ(n ; S)$ is distance magic.
\end{lemma}

\begin{proof}
Denote $\G = \Circ(n ; S)$ and set
\begin{equation}
\label{eq:def_la_mu_1}
	\lambda = b + d + 2n_0\quad \text{and}\quad \mu = c + d + 2n_0,
\end{equation}
where $\lambda$ and $\mu$ are treated as elements of $\ZZ_n$. Since $d$ is odd,~\eqref{eq:T1sys2} implies that $b+c+d$ is odd and is divisible by each of $d$, $d'$ and $d''$, and so $2(b+c+d) = 2n_0$ in $\ZZ_n$. It follows that 
\begin{equation}
\label{eq:sum1}
\lambda + 2\mu = d-b\quad \text{and}\quad 2\lambda + \mu = d-c\quad
\end{equation}
holds in $\ZZ_n$. Since $n_0$ is odd,~\eqref{eq:T1sys2} implies that each of $\lambda$ and $\mu$ is divisible by $4$. Moreover, since $\gcd(b,d) = 1$ and $d$ is coprime to both $d'$ and $d''$, $\gcd(\lambda,n) = 4d''$ and $\gcd(\mu,n) = 4d'$. 

Let $H = \la 4 \ra$ be the subgroup of $\ZZ_n$ generated by $4$ and let $\zeta \colon H \to \{0,1,\ldots, dd'-1\}$ and $\xi \colon H \to \{0,1,\ldots , d''-1\}$ be the functions such that
$$
	x = \zeta(x)\lambda + \xi(x)\mu 
$$
holds in $\ZZ_n$ for each $x \in H$. Since $\gcd(\lambda,n) = 4d''$ and $\gcd(\mu,n) = 4d'$ these two functions are indeed well defined. We now define a labeling $\ell_H$ of the elements of $H$ by setting 
$$
	\ell_H(x) = 1 + \zeta(x) + \xi(x)dd'\quad \text{for}\ x \in H.
$$ 
Observe that for each $x \in H$ we clearly have that $1 \leq \ell_H(x) \leq n_0$. Moreover, $\ell_H$ maps $H$ bijectively onto $\{1,2,\ldots , n_0\}$. 

Now, let $\delta \in \{-1,1\}$ be such that $n_0 \equiv \delta \pmod{4}$. We define a labeling $\ell$ on $\ZZ_{n}$ by the following rule:
\begin{equation}
\label{eq:T1_def_ell_2}
	\ell (x) = \left\{\begin{array}{ccc}
		\ell_H(x) & : & x \in H,\\
		n_0 + \ell_H(x-\delta n_0)    & : & x \in H + 1,\\
		4n_0+1-\ell_H(x+2n_0) & : & x \in H + 2,\\
		3n_0+1-\ell_H(x+\delta n_0) & : & x \in H + 3.\end{array}\right.
\end{equation}
It is clear that $\ell$ is a bijection from $\ZZ_{n}$ to $\{1,2,\ldots , n\}$. To see that it is in fact a distance magic labeling for $\G$ note first that $\ell(x) + \ell(x+2n_0) = n+1$ for each $x \in \ZZ_n$. Therefore, $w(x) + w(x+2n_0) = 6(n+1)$ holds for each $x \in \ZZ_n$, where $w(y)$ denotes the sum $\sum_{z \in N(y)}\ell(z)$. It thus suffices to prove that $w(x) = 3(n+1)$ for all $x \in \ZZ_n$ such that $x-d \in H \cup (H+1)$. 

To verify this, note that~\eqref{eq:def_la_mu_1} and~\eqref{eq:sum1} imply that for each $x \in \ZZ_n$ 
\begin{equation}
\label{eq:T1_2_expr}
\begin{array}{ccc}
	x-b & = & x-d + \lambda + 2\mu\\
	x-c & = & x-d + 2\lambda + \mu\\	
	x+b+2n_0 & = & x-d + \lambda\\
	x+c+2n_0 & = & x-d + \mu\\
	x+d+2n_0 & = & x-d + 2\lambda + 2\mu.
\end{array}
\end{equation}
Now, suppose $x \in \ZZ_n$ is such that $x-d \in H$. By~\eqref{eq:T1sys2} we have that $x-d, x-b, x-c \in H$ and $x+d, x+b, x+c \in H+2$, and so~\eqref{eq:T1_def_ell_2} implies that $w(x) = 3(n+1)$ if and only if 
\begin{equation}
\label{eq:T1_last}
\ell_H(x-d) + \ell_H(x-b) + \ell_H(x-c) = \ell_H(x+d+2n_0) + \ell_H(x+b+2n_0) + \ell_H(x+c+2n_0).
\end{equation}
Note that~\eqref{eq:T1_2_expr} implies that for each $y \in \{x-d, x-b, x-c\}$ there is precisely one $z \in \{x+d+2n_0, x+b+2n_0, x+c+2n_0\}$ such that $\zeta(y) = \zeta(z)$ and there is precisely one $z' \in \{x+d+2n_0, x+b+2n_0, x+c+2n_0\}$ such that $\xi(y) = \xi(z')$. This clearly shows that~\eqref{eq:T1_last} does indeed hold for all $x \in \ZZ_n$ with $x-d \in H$. The case when $x-d \in H+1$ is settled in a completely analogous way, where we simply substitute $x-\delta n_0$ for $x$ in~\eqref{eq:T1_2_expr} and then proceed as before.
\end{proof}

Combining together the results of this section we have the following classification of the distance magic circulants of valency $6$ for which all admissible characters are of type~1.

\begin{theorem}
\label{the:type1class}
Let $n \geq 7$ be an integer, let $S = \{\pm a, \pm b, \pm c\} \subset \ZZ_n$ be such that $|S| = 6$ and $\la S \ra = \ZZ_n$, and let $\G = \Circ(n; S)$. If all $j \in \A_n(S)$ are of type~1 then $\G$ is distance magic if and only if it is one of the graphs from Proposition~\ref{pro:type1}.
\end{theorem}

\section{Type~2}
\label{sec:type2}

In this section we consider admissible $j \in \A_n(S)$ of type~2 and classify the distance magic circulants of valency $6$ for which all admissible characters are of type~2.

\begin{lemma}
\label{le:type2}
Let $n \geq 7$ be an integer and let $S = \{\pm a, \pm b, \pm c\} \subset \ZZ_n$ be such that $|S| = 6$ and $\la S \ra = \ZZ_n$. Suppose at least one $j \in \A_n(S)$ of type~2 exists. Then $n = 3n_0$ for an integer $n_0 \geq 3$. Moreover, for each such $j$ there are $s_1, s_2, s_3 \in S$ with $\{\pm s_1, \pm s_2, \pm s_3 \} = S$ such that 
\begin{equation}
\label{eq:type2}
	j(s_2 - s_1) = n_0(1+3k_1) \quad \text{and}\quad j(s_3 - s_1) = n_0(2 + 3k_2)
\end{equation}
for some integers $k_1, k_2$, where these two equations are to be read within the ring of integers $\ZZ$.
\end{lemma}

\begin{proof}
Let $j \in \A_n(S)$ be of type~2. Observe that this implies that (modulo multiplying some of $a$, $b$ and $c$ by $-1$ and adding or subtracting multiples of $2\pi$) the angles $2\pi ja/n$, $2\pi jb/n$ and $2\pi jc/n$ are $2\pi/3$ apart from one another. In other words, for a suitable choice of $s_1 \in \{\pm a\}$, $s_2 \in \{\pm b\}$ and $s_3 \in \{\pm c\}$ we have that
$$
	\frac{2\pi js_2}{n} = \frac{2\pi js_1}{n} + \frac{2\pi}{3} + 2k_1\pi \quad \text{and}\quad \frac{2\pi js_3}{n} = \frac{2\pi js_1}{n} + \frac{4\pi}{3} + 2k_2\pi	
$$
for some integers $k_1$ and $k_2$. Rearranging we obtain
$$
	3j(s_2-s_1) = n(1+3k_1)\quad \text{and}\quad 3j(s_3-s_1) = n(2 + 3k_2),
$$
and so $n = 3n_0$ for some $n_0 \geq 3$. The above two equations then transform into~\eqref{eq:type2}.
\end{proof}

We can now characterize the only candidates for distance magic $6$-valent circulants for which all admissible characters are of type~2. 

\begin{proposition}
\label{pro:type2}
Let $n \geq 7$ be an integer and let $S = \{\pm a, \pm b, \pm c\} \subset \ZZ_n$ be such that $|S| = 6$ and $\la S \ra = \ZZ_n$. If $\G = \Circ(n; S)$ is distance magic and all $j \in \A_n(S)$ are of type~2, then $n = 3n_0$ for some $n_0 \geq 3$. Moreover, either $\G$ is trivial in the sense that there exists some $q \in \ZZ_n^*$ such that $qS = \{\pm 1, \pm (n_0 - 1), \pm (n_0 + 1)\}$ in which case $\G \cong C_{n_0}[3K_1]$, or the following both hold:
\begin{itemize}
\item $n_0 = dd'$ for coprime $d$ and $d'$ with $1 < d < d'$ both of which are coprime to $3$;
\item letting $\delta \in \{-1, 1\}$ be such that $n_0 \equiv \delta \pmod{3}$ and letting $c' \in \{1,2,\ldots , n-1\}$ be the unique solution of the system of congruences
\begin{equation}
\label{eq:type2system}
	\begin{array}{r@{\,}c@{\,}r@{}l}
	c' & \equiv & 0 & \pmod{3}\\
	c' & \equiv & 1 & \pmod{d}\\
	c' & \equiv & -1 & \pmod{d'},
	\end{array}
\end{equation}
there exists a $q \in \ZZ_n^*$ such that $qS = \{\pm 1, \pm (n_0 + \delta), \pm c'\}$.
\end{itemize}
\end{proposition}

\begin{proof}
Suppose that $\G$ is distance magic and all $j \in \A_n(S)$ (abbreviated to $\A_n$ in the rest of this proof) are of type~2. Lemma~\ref{le:type2} thus implies that $n = 3n_0$ for some $n_0 \geq 3$, while Lemma~\ref{le:common_divisor} implies that there is no divisor $d > 1$ of $n$ dividing all $j \in \A_n$. Suppose in addition that $\G$ is not the trivial example $C_{n_0}[3K_1]$. As in the proof of Proposition~\ref{pro:type1} we proceed by proving a series of claims. We also adopt the agreement that for $j, j' \in \A_n$ we let $s_1, s'_1, s_2, s'_2, s_3, s'_3 \in S$ and $k_1, k'_1, k_2, k'_2 \in \ZZ$ be as in~\eqref{eq:type2}.\smallskip

\noindent
{\sc Claim 1:} For any prime divisor $p$ of $n_0$ there are $\delta_p, \delta'_p \in \{-1, 1\}$ such that the whole $p$-part of $n_0$ divides each of $b - \delta_p a$, $c - \delta'_p a$ and $b - \delta_p\delta'_p c$. \\
Let $p$ be a prime divisor of $n_0$ and let $p^t$ be the $p$-part of $n_0$. By Lemma~\ref{le:common_divisor} there is $j \in \A_n$, which is coprime to $p$. Then Lemma~\ref{le:type2} implies that $p^t$ divides $s_2-s_1$ and $s_3 - s_1$, and consequently also $s_3-s_2$. The claim thus follows. 
\smallskip

\noindent 
{\sc Claim 2:} For any odd prime divisor $p$ of $n_0$ and for any $s, s' \in S$ with $s' \neq \pm s$, precisely one of $s \pm s'$ is divisible by $p$.\\
This follows from Claim~1 and the fact that if $p$ divided both $s + s'$ and $s - s'$, it would have to divide $s$ and $s'$ (as $p$ is odd), and so by Claim~1 it would divide each of $a$, $b$ and $c$, contradicting the assumption that $S$ generates $\ZZ_n$.
\smallskip

\noindent
{\sc Claim 3:} $n_0$ is not divisible by $3$.\\
By way of contradiction suppose $3$ divides $n_0$. By Claim~1 there are $\delta_3, \delta'_3 \in \{-1, 1\}$ such that the $3$-part of $n_0$ divides each of $b - \delta_3 a$, $c - \delta'_3 a$ and $b - \delta_3\delta'_3 c$. Then $a, b$ and $c$ are all coprime to $3$ (since $S$ generates $\ZZ_n$), and so $a$, $\delta_3 b$ and $\delta'_3 c$ all have the same nonzero remainder modulo $3$. Now, suppose there exists a $j \in \A_n$ such that $s_2 - s_1$ is not one of $\pm(b - \delta_3 a), \pm(c - \delta'_3 a)$ and $\pm(b - \delta_3\delta'_3 c)$. By Claim~2 we then have that $s_2 - s_1$ is coprime to $3$, and so~\eqref{eq:type2} implies that the $3$-part of $j$ equals the $3$-part of $n_0$. Therefore, $s_3 - s_1$ is also coprime to $3$. Since none of $s_1, s_2$ and $s_3$ is divisible by $3$, this implies that $s_2 \equiv s_3 \pmod{3}$ and hence also $s_2-s_1 \equiv s_3-s_1 \pmod{3}$. But then $(s_2-s_1) + (s_3-s_1)$ is also not divisible by $3$, which is impossible since~\eqref{eq:type2} yields $j(s_2-2s_1+s_3) = 3n_0(1+k_1+k_2)$ (recall that the $3$-part of $j$ equals the $3$-part of $n_0$). This (and an analogous proof for $s_3-s_1$) finally shows that for each $j \in \A_n$ each of $s_2-s_1$ and $s_3 - s_1$ is one of $\pm(b - \delta_3 a), \pm(c - \delta'_3 a)$ and $\pm(b - \delta_3\delta'_3 c)$. Combining this with Lemma~\ref{le:common_divisor} and Claim~2 we thus find that the whole $n_0$ divides each of $b - \delta_3 a$, $c - \delta'_3 a$ and $b - \delta_3\delta'_3 c$. It thus follows that $S = \{\pm a, \pm (a + n_0), \pm (a + 2n_0)\}$. But since $3$ divides $n_0$ and $S$ generates $\ZZ_n$, it must be that $\gcd(a,n) = 1$, and so multiplication by a suitable $q \in \ZZ_n^*$ gives $qS = \{\pm 1, \pm (n_0+1), \pm (2n_0+1)\} = \{\pm 1, \pm (n_0-1), \pm (n_0+1)\}$, contradicting the assumption that $\G$ is not trivial.
\smallskip

\noindent
{\sc Claim 4:} For precisely one $s \in \{a,b,c\}$ we have that $\gcd(n,s) = 3$, while the remaining two are coprime to $n$.\\ 
That no prime $p$, other than perhaps $3$, can divide any of $a$, $b$ and $c$ follows from Claim~1 and the fact that $S$ generates $\ZZ_n$. Now, take any $j \in \A_n$. By Claim~3 and~\eqref{eq:type2} none of $j$, $s_2-s_1$ and $s_3-s_1$ is divisible by $3$. Moreover, as $j(s_3-s_2) = n_0(1 + 3k_2 - 3k_1)$, $s_3-s_2$ is also not divisible by $3$. But then precisely one of $s_1$, $s_2$ and $s_3$ must be divisible by $3$. 
\smallskip

We now complete the proof as follows. Claim~4 implies that multiplying by a suitable $q \in \ZZ_n^*$ we can assume that $S = \{\pm 1, \pm b, \pm c\}$, where one of $b$ and $c$ is coprime to $n$, while the other is divisible by $3$ but coprime to $n_0$. With no loss of generality assume $b$ is coprime to $n$. Then precisely one of $b \pm 1$ is coprime to $3$. Exchanging the roles of $b$ and $-b$ if necessary we can thus assume that $b + 1$ is coprime to $3$. By Lemma~\ref{le:type2} it follows that for any $j \in \A_n$ one of $s_2-s_1$, $s_3-s_1$ and $s_2-s_3$ (which are all coprime to $3$ by Claim~3 and \eqref{eq:type2}) must be one of $\pm(b + 1)$. Then Lemma~\ref{le:common_divisor} and~\eqref{eq:type2} imply that for each prime divisor $p$ of $n_0$ the whole $p$-part of $n_0$ divides $b + 1$, and so $n_0$ divides $b + 1$. As $b \neq \pm 1$, we must have that $b$ is one of $-1 \pm n_0$, depending on which of the two is coprime to $3$. In other words, letting $\delta \in \{-1,1\}$ be such that $n_0 \equiv \delta \pmod{3}$ we have that $b = \pm (n_0 + \delta)$. Therefore, $S = \{\pm 1, \pm (n_0 + \delta), \pm c\}$ for some $c$ divisible by $3$ but coprime to $n_0$.

Claim~1 implies that for any prime divisor $p$ of $n_0$ the whole $p$-part of $n_0$ divides one of $c \pm 1$ (in the case that $p = 2$ and the $2$-part of $n_0$ is $2$, it divides both). If for all prime divisors $p$ of $n_0$ the $p$-part of $n_0$ divides $c - 1$ (or analogously if it always divides $c + 1$), then $c = \pm (n_0 - \delta)$, contradicting the assumption that $\G$ is not trivial. We now define coprime $d$ and $d'$ with $n_0 = dd'$. Let $p$ be any prime divisor of $n_0$. If $p \neq 2$, then by Claims~1 and 2 the $p$-part of $n_0$, say $p^t$, divides precisely one of $c-1$ and $c+1$. We let $p^t$ be the $p$-part of $d$ or $d'$, depending on whether $p^t$ divides $c-1$ or $c+1$, respectively. If $p = 2$, we do the following. If the $p$-part of $n_0$ is $2$, then we let the $2$-part of $d$ be $2$ (and $d'$ be odd). If however the $2$-part $2^t$ of $n_0$ is at least $4$, then we let $2^t$ be the $2$-part of $d$ or $d'$, depending on whether $2^t$ divides $c-1$ or $c+1$, respectively. Clearly, $n_0 = dd'$. Moreover, we have that $3 \mid c$, $c \equiv 1 \pmod{d}$ and $c \equiv -1 \pmod{d'}$ which determines $c$ completely. Exchanging the roles of $c$ and $-c$ if necessary we can assume that $d < d'$. 
\end{proof}

We next show that in fact each distance magic circulant $\Circ(n ; \{\pm a, \pm b, \pm c\})$ of valency $6$ for which all admissible characters are of type~2 is of odd order.

\begin{corollary}
\label{cor:type2}
Let $n \geq 7$ be an integer and let $S = \{\pm a, \pm b, \pm c\} \subset \ZZ_n$ be such that $|S| = 6$ and $\la S \ra = \ZZ_n$. If $\G = \Circ(n; S)$ is distance magic and all $j \in \A_n(S)$ are of type~2, then $n$ is odd.
\end{corollary}

\begin{proof}
In view of Proposition~\ref{pro:type2} we know that $n = 3n_0$ for some $n_0 \geq 3$. By way of contradiction suppose $n_0$ is even. We first deal with the trivial example $\Circ(n ; \{\pm 1, \pm (n_0 - 1), \pm (n_0 + 1)\})$. Since $\G$ is distance magic, we cannot have $n_0 \equiv 2 \pmod{4}$. This is easy to see (but see also~\cite{ShaAliSim09}). Namely, since $n$ is even the magic constant $3(n+1)$ of $\G$ is odd. But since the set of neighbors of a vertex $x$ consists of the union of cosets $H + x + 1$ and $H + x - 1$ of the subgroup $H = \la n_0 \ra \leq \ZZ_n$, the fact that the number of cosets of $H$ in $\ZZ_n$ is twice an odd number shows that the sum of the labels on each coset of $H$ would have to be $3(n+1)/2$, which is not an integer. If $n_0$ is divisible by $4$, say $n_0 = 4m$ for some $m \geq 1$, then $\G$ is in fact indeed distance magic (see~\cite{ShaAliSim09}) but in this case not all admissible $j \in \A_n$ are of type $2$. Namely, since $j = 3m = n/4$ satisfies~\eqref{eq:cos_sum}, it is admissible. However, if $j$ satisfied \eqref{eq:type2} for some $s_1, s_2, s_3 \in S$ and integers $k_1, k_2$, then $3(s_2-s_1) = 4(1+3k_1)$ would hold, which is clearly impossible. Therefore, this $j$ is not of type~2.

To complete the proof we now assume that there is no $q \in \ZZ_n ^*$ such that $qS$ is $\{\pm 1, \pm (n_0 - 1), \pm (n_0 + 1)\}$. By Claim~2 from the proof of Proposition~\ref{pro:type2} we then have that $n_0$ is coprime to $3$. Since all $j \in \A_n(S)$ are of type~2, we see that~\eqref{eq:type2} implies that all $j \in \A_n(S)$ are coprime to $3$. But then~\eqref{eq:characters} implies that for each $j \in \A_n(S)$ we get
\begin{equation}
\label{eq:sub_ind_3}
	\sum_{x \in \la n_0 \ra}\chi_j(x) = 1 + e^{\frac{2 \pi j\mathbf{i}}{3}} + e^{\frac{4 \pi j\mathbf{i}}{3}} = 1 + e^{\frac{2 \pi \mathbf{i}}{3}} + e^{\frac{4 \pi \mathbf{i}}{3}} = 0.
\end{equation}
The discussion from the last paragraph of Section~\ref{sec:Prelim} thus shows that for any eigenvector for the eigenvalue $0$ of the adjacency matrix of $\G$ the sum of the entries corresponding to $0$, $n_0$ and $2n_0$ is $0$. But this contradicts Lemma~\ref{le:eig0} since the nominators of the numbers in~\eqref{eq:seq} are all odd (recall that $n$ was assumed to be even). This contradiction shows that either $\G$ is not distance magic or there are admissible $j \in \A_n(S)$ which are not of type~2.
\end{proof}

We next show that the examples from the second part of Proposition~\ref{pro:type2} with $n_0$ odd are all distance magic.

\begin{lemma}
\label{le:T2_labeling}
Let $d$ and $d'$, where $1 < d < d'$, be coprime integers which are both coprime to $6$, and let $n_0 = dd'$ and $n = 3n_0$. Let $\delta \in \{-1, 1\}$ be such that $n_0 \equiv \delta \pmod{3}$, let $c \in \{1,2,\ldots , n-1\}$ be the unique solution of~\eqref{eq:type2system} and set $S = \{\pm 1, \pm (n_0 + \delta), \pm c\}$. Then the circulant $\Circ(n ; S)$ is distance magic.
\end{lemma}

\begin{proof}
Denote $\G = \Circ(n ; S)$ and set
\begin{equation}
\label{eq:def_la_mu_2}
	\lambda = c + 1 -\delta n_0\quad \text{and}\quad \mu = c - 1 + \delta n_0,
\end{equation}
where $\lambda$ and $\mu$ are treated as elements of $\ZZ_n$. Note that $\delta n_0 \equiv 1 \pmod{3}$, and so as $d$ and $d'$ are odd and coprime,~\eqref{eq:type2system} implies that 
$$
	\gcd(\lambda, n) = 3d' \quad \text{and} \quad \gcd(\mu, n) = 3d.
$$
Similarly as in the proof of Lemma~\ref{le:T1DMcase2} let $H = \la 3 \ra \leq \ZZ_n$ and let $\zeta \colon H \to \{0,1,\ldots, d-1\}$ and $\xi \colon H \to \{0,1,\ldots , d'-1\}$ be the functions such that
$$
	x = \zeta(x)\lambda + \xi(x)\mu 
$$
holds in $\ZZ_n$ for each $x \in H$. We define a labeling $\ell_H$ of the elements of $H$ by setting 
\begin{equation}
\label{eq:ell_H_T2}
	\ell_H(x) = 1 + \zeta(x) + \xi(x)d\quad \text{for}\ x \in H.
\end{equation}
Observe that for each $x \in H$ we have that $1 \leq \ell_H(x) \leq n_0$. Moreover, $\ell_H$ maps $H$ bijectively onto $\{1,2,\ldots , n_0\}$. We now define a labeling $\ell \colon \ZZ_n \to \{1,2,\ldots , n\}$ of the vertices of $\G$ as follows:
\begin{equation}
\label{eq:T2_def_ell}
	\ell (x) = \left\{\begin{array}{ccc}
		\ell_H(x) & : & x \in H,\\
		n_0 + \ell_H(x+\lambda-1)    & : & x \in H + 1,\\
		2n_0+\ell_H(-2(x+2(\lambda-1))) & : & x \in H + 2.\end{array}\right.
\end{equation}
Since $3$ divides $\lambda$, we have that $\lambda-1 \equiv -1 \pmod{3}$, and so $\ell$ is well defined. Moreover, as $n$ is odd, $-2$ has a multiplicative inverse in $\ZZ_n$, and so $\ell$ is a bijection. We prove that it is in fact a distance magic labeling by proving a series of claims.
\smallskip

\noindent
{\sc Claim 1:} For each $q \in \ZZ_n ^*$ and for each $x \in H$ we have that 
$$	
	\ell_H(q(x-c)) + \ell_H(q(x+c)) = \ell_H(q(x+\delta n_0-1)) + \ell_H(q(x-\delta n_0 + 1)). 
$$
Note first that the claim makes sense since $c, \delta n_0-1 \in H$, so that we can actually compute $\ell_H$ on all of these elements. To see that the equality does indeed hold it suffices to observe that by~\eqref{eq:def_la_mu_2} each of the following holds in $\ZZ_n$:
$$
\begin{array}{ccc}
	x+c & = & x-c + \lambda + \mu\\
	x+\delta n_0-1 & = & x-c + \mu\\	
	x-\delta n_0 + 1 & = & x-c + \lambda.
\end{array}
$$
\smallskip

\noindent
{\sc Claim 2:} For each $x \in \ZZ_n$ we have that $\ell(x-c) + \ell(x+c) = \ell(x+\delta n_0 - 1) + \ell(x-\delta n_0 + 1)$.\\
This is an immediate consequence of Claim~1 and \eqref{eq:T2_def_ell}.
\smallskip

\noindent
{\sc Claim 3:} We have  $\zeta(c) = (d+1)/2$ and $\xi(c) = (d'+1)/2$.\\
This holds as $n$ is odd, $\lambda + \mu \equiv 2c \pmod{n}$ and $d\lambda \equiv 0 \pmod{n}$ and $d'\mu \equiv 0 \pmod{n}$. 
\smallskip

\noindent
{\sc Claim 4:} If for each integer $x$ we let $(x)_{d} \in \{0,1,\ldots , d-1\}$ be the remainder of $x$ modulo $d$, then for each $y \in \ZZ_n$ we have that $(y)_d + (y+(d+1)/2)_d + (-2y-2)_d = d + (d-3)/2$.\\
Denote the left-hand sum by $z$. Observe first that $z \equiv (d-3)/2 \pmod{d}$. It thus suffices to verify that $(d-3)/2 < z < 2d+(d-3)/2$. Since $d \geq 5$ and is odd we clearly have that 
$$
(d-3)/2 < (y)_d + (y+(d+1)/2)_d \leq d + (d-3)/2,
$$
and so the claim follows.
\smallskip

\noindent
{\sc Claim 5:} For any $x \in \ZZ_n$ we have that $\ell(x) + \ell(x+\delta n_0) + \ell(x-\delta n_0) = 3(n+1)/2$.\\
Observe first that (since $\delta n_0 \equiv 1 \pmod{3}$) it suffices to prove that this holds for each $x \in H$. Suppose then that $x \in H$. By~\eqref{eq:def_la_mu_2} and \eqref{eq:T2_def_ell} and since $\delta n_0 \equiv 1 \pmod{3}$ we have that 
$$
	\ell(x+\delta n_0) = n_0 + \ell_H(x+\delta n_0 + \lambda - 1) = n_0 + \ell_H(x+c).
$$
Similarly, since $3n_0 = 0$ in $\ZZ_n$, we have that
$$
	\ell(x-\delta n_0) = 2n_0 + \ell_H(-2(x-\delta n_0 + 2(\lambda - 1))) = 2n_0 + \ell_H(-2x-4c).
$$
By \eqref{eq:ell_H_T2} we thus have that $\ell(x) + \ell(x+\delta n_0) + \ell(x-\delta n_0)$ equals
$$
	 3+n+\zeta(x) + \zeta(x+c) + \zeta(-2x-4c) + d(\xi(x) + \xi(x+c) + \xi(-2x-4c)).
$$
To prove our claim it thus suffices to show that this equals 
$$
	3+n+d + (d-3)/2 + d(d' + (d'-3)/2) = 3(n+1)/2.
$$
This for sure holds if we can prove that
$$
	\zeta(x) + \zeta(x+c) + \zeta(-2x-4c) = d + (d-3)/2 \quad \text{and}\quad \xi(x) + \xi(x+c) + \xi(-2x-4c) = d' + (d'-3)/2.
$$
By Claim~3 we have that 
$$
	\zeta(x+c) = (\zeta(x) + (d+1)/2)_d \quad \text{and} \quad \zeta(-2x-4c) = (-2\zeta(x)-4(d+1)/2)_d = (-2\zeta(x)-2)_d.
$$
By Claim~4 it thus follows that 
$$
	\zeta(x) + \zeta(x+c) + \zeta(-2x-4c) = d + (d-3)/2,
$$
as claimed. A similar argument in which we replace $d$ by $d'$ shows that also $\xi(x) + \xi(x+c) + \xi(-2x-4c) = d' + (d'-3)/2$ holds.
\smallskip

\noindent
To complete the proof let $x \in \ZZ_n$ and observe that $\pm(n_0+\delta) = \pm(\delta n_0 + 1)$. By Claim~2 we have that the sum of the labels of the neighbors of $x$ is 
$$
	\ell(x-1) + \ell(x+1) + \ell(x+\delta n_0 + 1) + \ell(x - \delta n_0 - 1) + \ell(x+\delta n_0 - 1) + \ell(x - \delta n_0 + 1),
$$
which by Claim~5 equals $3(n+1)/2 + 3(n+1)/2 = 3(n+1)$, as required.
\end{proof}

Note that, in view of Proposition~\ref{pro:type2}, Corollary~\ref{cor:type2} and Lemma~\ref{le:T2_labeling} the only examples with all admissible characters being of type~2 for which we have not yet given a distance magic labeling are the trivial examples $\Circ(3n_0 ; \{\pm 1, \pm (n_0 - 1), \pm (n_0 + 1)\}) \cong C_{n_0}[3K_1]$ where $n_0$ is odd. Recall from Section~\ref{sec:Prelim} that these graphs are indeed distance magic (see~\cite{ShaAliSim09}). It is interesting to note that, at least when $n_0 = dd'$, where $1 < d < d'$ are coprime and are both also coprime to $6$, a distance magic labeling can be obtained by simply taking the labeling from the proof of Lemma~\ref{le:T2_labeling} (this follows from Claim~2 of that proof). 

We wrap up this section by stating the obtained classification of the distance magic circulants of valency $6$ for which all admissible characters are of type~2.

\begin{theorem}
\label{the:type2class}
Let $n \geq 7$ be an integer, let $S = \{\pm a, \pm b, \pm c\} \subset \ZZ_n$ be such that $|S| = 6$ and $\la S \ra = \ZZ_n$, and let $\G = \Circ(n; S)$. If all $j \in \A_n(S)$ are of type~2 then $\G$ is distance magic if and only if $n$ is odd and $\G$ is one of the graphs from Proposition~\ref{pro:type2}.
\end{theorem}

\section{Type~3}
\label{sec:type3}

In this section we analyze the examples which admit at least one exceptional solution of~\eqref{eq:Cox}. In other words, we analyze the examples for which at least one admissible character is of type~3.

\begin{lemma}
\label{le:type3}
Let $n \geq 7$ be an integer and let $S = \{\pm a, \pm b, \pm c\} \subset \ZZ_n$ be such that $|S| = 6$ and $\la S \ra = \ZZ_n$. Suppose there exists a $j \in \A_n(S)$ of type~3. Then the following all hold:
\begin{itemize}
\item $n = 30n_0$ for an integer $n_0 \geq 1$;
\item $j = j_0n_0$ for some $j_0 \in \{1,7,11,13,17,19,23,29\}$;
\item up to multiplying by $-1$ and changing the roles of $a$, $b$ and $c$ one of the following holds:
   \begin{itemize}
   \item $a \equiv 3 \pmod{30}$, $b \equiv 9 \pmod{30}$ and $c \equiv 10 \pmod{30}$, or
   \item $a \equiv 5 \pmod{30}$, $b \equiv 6 \pmod{30}$ and $c \equiv 12 \pmod{30}$.
   \end{itemize}
\end{itemize}
\end{lemma}

\begin{proof}
Let $j \in \A_n(S)$ be of type~3 and suppose that the corresponding solution of~\eqref{eq:Cox} is of the first of the two forms from~\eqref{eq:sol3}. With no loss of generality we can then assume that
$$
	\frac{2\pi ja}{n} = \frac{\pi}{5} + 2k_1\pi,\quad \frac{2\pi jb}{n} = \frac{3\pi}{5} + 2k_2\pi\quad \text{and}\quad \frac{2\pi jc}{n} = \frac{2\pi}{3} + 2k_3\pi	
$$
for some $k_1, k_2, k_3 \in \ZZ$. Rearranging we obtain
\begin{equation}
\label{eq:type3aux}
	10ja = n(1+10k_1), \quad 10jb = n(3+10k_2)\quad \text{and}\quad 3jc = n(1+3k_3).
\end{equation}
The first of these shows that $10 \mid n$ and the third that $3 \mid n$, and so $n = 30n_0$ for some $n_0 \geq 1$. Using this in~\eqref{eq:type3aux}, dividing the first two equations by 10 and the third by 3, we find that $n_0$ divides each of $ja$, $jb$ and $jc$. Since $\la S \ra = \ZZ_n$ it follows that $n_0 \mid j$, and so $j = j_0n_0$ for some positive integer $j_0$. This translates~\eqref{eq:type3aux} to
$$
	j_0a = 3(1+10k_1),\quad j_0b = 3(3+10k_2)\quad \text{and}\quad j_0c = 10(1+3k_3).
$$
The third of these implies that $3$ divides none of $j_0$ and $c$, and so the first two equations imply that both $a$ and $b$ are divisible by $3$. Similarly, $a$, $b$ and $j_0$ are all coprime to $10$, and consequently $c$ is divisible by $10$. Observe also that as $j < n$ and $\gcd(j_0, 30) = 1$ this in fact implies that $j_0 \in \{1,7,11,13,17,19,23,29\}$, as claimed. It is now clear that $c \equiv \pm 10 \pmod{30}$ and that, after possibly replacing the roles of $a$ and $b$, we have that $a \equiv \pm 3 \pmod{30}$ and $b \equiv \pm 9 \pmod{30}$.

The proof for the solution of~\eqref{eq:Cox} of the second of the two forms from~\eqref{eq:sol3} is done in a completely analogous way and is left to the reader.
\end{proof}

\begin{proposition}
\label{pro:type3}
Let $n \geq 7$ be an integer and let $S = \{\pm a, \pm b, \pm c\} \subset \ZZ_n$ be such that $|S| = 6$ and $\la S \ra = \ZZ_n$. Suppose that the circulant $\G = \Circ(n; S)$ is distance magic and that there exists at least one $j \in \A_n(S)$ of type~3. Then the following all hold:
\begin{itemize}
\item $n$ is divisible by $60$, the $3$-part of $n$ is $3$ and the $5$-part of $n$ is $5$;
\item there is a unique $s \in S$ with $5 \mid s$ and $s < n/2$, and this $s$ is one of $n/12$ and $5n/12$.
\item there exists at least one $j \in \A_n(S)$ of type~1 and there are no $j \in \A_n(S)$ of type~2;
\item all $j \in \A_n(S)$ of type~1 are odd and are divisible by $15$, while all $j \in \A_n(S)$ of type~3 are even and coprime to $15$;
\end{itemize}
\end{proposition}

\begin{proof}
By Lemma~\ref{le:type3} precisely one of $a$, $b$, $c$ is coprime to $3$, precisely one is divisible by $5$, and $n = 30n_0$ for some $n_0 \geq 1$. We abbreviate $\A_n(S)$ to $\A_n$ and proceed by a series of claims.
\smallskip

\noindent
{\sc Claim~1:} Each $j \in \A_n$ of type~1 is divisible by $3$.\\
Suppose $j \in \A_n$ is of type~1. Letting $s_1, s_2, s_3 \in S$ be as in Lemma~\ref{le:type1} we then have that $js_2$ and $j(s_1+s_3)$ are both divisible by $3$. Since precisely one of $s_1, s_2$ and $s_3$ is coprime to $3$, it thus follows that $3 \mid j$, as claimed. 
\smallskip

\noindent
{\sc Claim~2:} $n_0$ is coprime to $3$.\\
By way of contradiction suppose $3 \mid n_0$. By Claim~1 all admissible $j \in \A_n$ of type~1 are divisible by 3, while by Lemma~\ref{le:type3} the same holds for each $j \in \A_n$ of type~3. By Lemma~\ref{le:common_divisor} there thus must exist at least one $j \in \A_n$ of type~2 which is coprime to $3$. Letting $s_1, s_2$ and $s_3$ be as in Lemma~\ref{le:type2} for this $j \in \A_n$ we have that $3$ divides each of $s_2-s_1$ and $s_3-s_1$, which contradicts the fact that precisely one of $s_1, s_2$ and $s_3$ is coprime to $3$.
\smallskip

\noindent
{\sc Claim~3:} There are no $j \in \A_n$ of type~2.\\
Suppose to the contrary that there exists a $j \in \A_n$ of type~2 and let $s_1$, $s_2$ and $s_3$ be as in Lemma~\ref{le:type2}. By Claim~2 and Lemma~\ref{le:type2} none of $s_2-s_1$, $s_3-s_1$ and $s_3-s_2$ is divisible by $3$, contradicting the fact that two of $s_1$, $s_2$ and $s_3$ are divisible by $3$. 
\smallskip

\noindent
{\sc Claim~4:} There exists at least one $j \in \A_n$ with $3 \mid j$.\\
By way of contradiction suppose that all $j \in \A_n$ are coprime to $3$. Letting $H = \la 10n_0\ra$ be the subgroup of $\ZZ_n$ of order $3$ we then see, just as in~\eqref{eq:sub_ind_3}, that $\sum_{x \in H}\chi_j(x) = 0$ holds for each $j \in \A_n$. As in the proof of Corollary~\ref{cor:type2} we see (as $n$ is even) that this contradicts Lemma~\ref{le:eig0}. 
\smallskip

\noindent
{\sc Claim~5:} There exists at least one $j \in \A_n$ of type~1 and $n_0$ is even.\\
Claim~2 implies that $n$ is not divisible by $9$, and so Lemma~\ref{le:type3} shows that all $j \in \A_n$ of type~3 are coprime to $3$. Claims~3 and 4 therefore force the existence of at least one $j \in \A_n$ of type~1. Lemma~\ref{le:type1} then yields $4 \mid n$, and consequently $n_0$ is even. 
\smallskip

\noindent
{\sc Claim~6:} Each $j \in \A_n$ of type~1 is divisible by $5$ and $n_0$ is coprime to $5$.\\
Take any $j \in \A_n$ of type~1 and let $s_1$, $s_2$ and $s_3$ be as in Lemma~\ref{le:type1}. Since $s_1 + s_3$ must be even, $s_1$ and $s_3$ are of the same parity, and so Lemma~\ref{le:type3} implies that $s_2$ is divisible by $5$ while $s_1$ and $s_3$ are not. In fact, the possibilities for $a$, $b$ and $c$ from Lemma~\ref{le:type3} show that $s_1+s_3$ is not divisible by $5$. Lemma~\ref{le:type1} thus implies that $j$ must be divisible by $5$, which proves the first part of this claim. But since we now know that each $j \in \A_n$ of type~1 is divisible by $5$, Lemma~\ref{le:common_divisor}, Lemma~\ref{le:type3} and Claim~3 imply that $n_0$ is coprime to $5$.
\smallskip

\noindent
{\sc Claim~7:} Each $j \in \A_n$ of type~1 is odd and the corresponding $s_2$ from Lemma~\ref{le:type1} is divisible by $5$.\\
Let $j \in \A_n$ be of type~1 and let $s_1, s_2$ and $s_3$ be as in Lemma~\ref{le:type1}. Since $s_1 + s_3$ is even, $s_1$ and $s_3$ are of the same parity, and so Lemma~\ref{le:type3} implies that $s_2$ is divisible by $5$. To prove the first part of the claim suppose that $j$ is even and denote $m = n/4$, which must be even by~\eqref{eq:type1}. Lemma~\ref{le:type3} then implies that all $j' \in \A_n$ of type~3 are even. By Lemma~\ref{le:common_divisor} and Claim~3 there thus exists an odd $j'' \in \A_n$ of type 1. The corresponding $s_2''$ from Lemma~\ref{le:type1} then must be even. But since $s_2''$ is divisible by $5$, Lemma~\ref{le:type3} implies that $s_2 = \pm s_2''$, and so $js_2$ and $j''s_2''$ do not have the same $2$-part (since $j$ is even while $j''$ is odd), contradicting~\eqref{eq:type1}. This thus shows that no even $j \in \A_n$ of type~1 exists. 
\smallskip

\noindent
To complete the proof let $s \in S$ be the unique element with $s < n/2$ and $5 \mid s$ (see Lemma~\ref{le:type3}). By Claim~7 it then follows that for each $j \in \A_n$ of type~1 the corresponding $s_2$ is $\pm s$. Now, let $p$ be any prime divisor of $n$, different from $2$, $3$ and $5$. By Lemma~\ref{le:type3} this $p$ divides each $j' \in \A_n$ of type~3, and so Lemma~\ref{le:common_divisor} and Claim~3 imply that there is at least one $j \in \A_n$ of type~1, which is coprime to $p$. By Lemma~\ref{le:type1} it thus follows that the whole $p$-part of $n$ divides $s$. Moreover, by Claim~7 the $2$-part of $s$ coincides with the $2$-part of $n/4$. Since $s$ is not divisible by $3$ (by Lemma~\ref{le:type3}) we finally see that $s$ is one of $n/12$ and $5n/12$. 
\end{proof}

Since Proposition~\ref{pro:type3} shows that there is no distance magic circulant of valency $6$ all of whose admissible characters are of type~3, the previous two sections provide a complete classification of distance magic circulants of valency $6$ for which all admissible characters are of the same type. In fact, combining together Lemmas~\ref{le:type1} and \ref{le:type2}, Theorems~\ref{the:type1class} and~\ref{the:type2class} and Proposition~\ref{pro:type3} we have the following result (which in particular implies that there are no connected distance magic circulants of valency $6$ and order twice an odd number).

\begin{theorem}
\label{the:class_part}
Let $n \geq 7$ be an integer, let $S = \{\pm a, \pm b, \pm c\} \subset \ZZ_n$ be such that $|S| = 6$ and $\la S \ra = \ZZ_n$, and let $\G = \Circ(n; S)$. Then the following both hold:
\begin{itemize}
\item If $n$ is not divisible by $3$ then $\G$ is distance magic if and only if $\G$ is one of the graphs from Proposition~\ref{pro:type1}.
\item If $n$ is not divisible by $4$ then $\G$ is distance magic if and only if $n$ is odd and $\G$ is one of the graphs from Proposition~\ref{pro:type2}.
\end{itemize}
\end{theorem}

\section{Mixed types}
\label{sec:mixed}

By the results of the previous sections the only remaining circulants of valency $6$ that need to be investigated are those having admissible characters of at least two types. By Proposition~\ref{pro:type3} all such distance magic examples have admissible characters of type~1, and so Lemmas~\ref{le:type1}, \ref{le:type2} and~\ref{le:type3} imply that for such graphs the order $n$ needs to be divisible by $12$. Moreover, either all admissible characters which are not of type~1 are of type 2 or they are all of type 3. 

Using Lemmas~\ref{le:type1} and~\ref{le:type2} it is easy to see that, up to isomorphism, the only two connected circulants of order $12$ and valency $6$ for which there are admissible characters of types 1 and 2 are the trivial examples $\Circ(12 ; \{\pm 1, \pm 3, \pm 5\}) \cong \mathrm{Ml}_{3}[2K_1]$ and $\Circ(12 ; \{\pm 2, \pm 3, \pm 4\})\cong \mathrm{Pr}_{3}[2K_1]$ which we already know are indeed distance magic.

The situation get more interesting already at the next possible order, namely $24$. Using a computer one can easily verify that up to isomorphism there are precisely five connected $6$-valent circulants of the form $\Circ(24; S)$ for which the set of all admissible $j \in \A_{24}(S)$ does not have a nontrivial common divisor (so as to present a possible distance magic graph in view of Lemma~\ref{le:common_divisor}). These are the trivial examples $\G_1 = \Circ(24; \{\pm 1, \pm 6, \pm 11\}) \cong \mathrm{Ml}_{6}[2K_1]$ and $\G_2 = \Circ(24; \{\pm 1, \pm 7, \pm 9\}) \cong C_8[3K_1]$ and the three nontrivial ones
$$
	\G_3 = \Circ(24;\{\pm 1, \pm 2, \pm 3\}),\  \G_4 = \Circ(24;\{\pm 1, \pm 3, \pm 10\}),\  \G_5 = \Circ(24;\{\pm 1, \pm 5, \pm 6\}).
$$
That $\G_1$ and $\G_2$ are distance magic was explained in Section~\ref{sec:type1}. As for the nontrivial examples, it is not difficult to verify that for each of the graphs $\G_3$, $\G_4$ and $\G_5$ the corresponding set $\A_{24}(S)$ is $\{\pm 3, \pm 8, \pm 9\}$, where $\pm 3$ and $\pm 9$ are of type 1, while $\pm 8$ are of type 2. Consequently, by the discussion at the end of Section~\ref{sec:Prelim} the eigenspace corresponding to the eigenvalue $0$ of the adjacency matrix of $\G_i$ is the same for all $i \in \{3,4,5\}$. Lemma \ref{le:eig0} therefore implies that if any of these graphs is distance magic then they all are (and any distance magic labeling for one of them is also a distance magic labeling for the remaining two). It follows from~\cite[Theorem 10]{Cic16} that $\G_3$ is distance magic. In fact, one can check that a magic labeling $\ell$ is given by
$$
\begin{array}{c c c c c c}
	\ell(0) = 2, & \ell(1) = 7, & \ell(2) = 15, & \ell(3) = 5, & \ell(4) = 22, & \ell(5) = 18,  \\
	\ell(6) = 11, & \ell(7) = 19 & \ell(8) = 3, & \ell(9) = 8, & \ell(10) = 13, & \ell(11) = 6, \\
	\ell(12) = 23, & \ell(13) = 16, & \ell(14) = 12, &  \ell(15) = 20, &  \ell(16) = 1, & \ell(17) = 9, \\
	\ell(18) = 14, & \ell(19) = 4, &  \ell(20) = 24, & \ell(21) = 17, & \ell(22) = 10, &  \ell(23) = 21.
\end{array}
$$
The graphs $\G_3$, $\G_4$ and $\G_5$ are thus all examples of nontrivial distance magic circulants of valency $6$ having admissible characters of types~1 and~2. It appears that the analysis of all distance magic circulants of valency $6$ having admissible characters of type~1 and also of type~2 is considerably more complicated than the analysis from Section~\ref{sec:type1} and~\ref{sec:type2}. We thus propose the following problem.

\begin{problem}
\label{prob1}
Classify distance magic circulants of valency $6$ for which some admissible characters are of type~1 and some are of type~2.
\end{problem}

Let us finally consider circulants of valency $6$ having admissible characters of types~1 and 3. By Proposition~\ref{pro:type3} the first possible order of such a circulant is $60$. That proposition also states that if $\G = \Circ(60; \{\pm a, \pm b, \pm c\})$ is such an example then we can assume that $a$ is one of $5$ and $25$. Lemma~\ref{le:type3} then further implies that, up to changing the roles of $b$ and $c$, we have that $b \in \{6, 24\}$ and $c \in \{12, 18\}$. Finally, by Lemma~\ref{le:type1} precisely one of $b$ and $c$ is divisible by $4$, and so there are only four possibilities for $a$, $b$ and $c$. It is easy to see (consider multiplication by $7$ and by $13$) that all four give isomorphic graphs, so we can just take one of them, say $\G = \Circ(60; \{\pm 5, \pm 6, \pm 12\})$. It is straightforward to verify that the corresponding $\A_{60}(S)$ is $\{ \pm 2, \pm 14, \pm 15, \pm 22, \pm 26\}$, where $j = \pm 15$ are of type 1, while all other $j \in \A_{60}(S)$ are of type 3. We do not know whether $\G$ is distance magic or not, and so we propose the following problem. 

\begin{problem}
Classify distance magic circulants of valency $6$ for which some admissible characters are of type~1 and some are of type~3.
\end{problem}

We conclude the paper by a comment that perhaps indicates that at least Problem~\ref{prob1} might be quite difficult. Using a computer one can verify that there are precisely $15$ nonisomorphic connected circulants of the form $\Circ(60; S)$ for which the set of all admissible $j \in \A_{60}(S)$ does not have a nontrivial common divisor (so as to make them potential candidates for being distance magic by Lemma~\ref{le:common_divisor}). Except for the above mentioned $\Circ(60; \{\pm 5, \pm 6, \pm 12\})$ all of the remaining $14$ have admissible $j \in \A_{60}(S)$ of type~1 and also admissible $j'\in \A_{60}(S)$ of type~2. Three of these are the trivial $\mathrm{Pr}_{15}[2K_1]$, $\mathrm{Ml}_{15}[2K_1]$ and $C_{20}[3K_1]$. But there are thus $11$ nontrivial examples. We did not try to see which of them are distance magic as this does not seem to be easy to determine. But the number of examples does explain our above mentioned feeling that the analysis of examples having admissible characters of types~1 and~2 might be difficult. We also mention that two of them, for instance $\Circ(60;\{\pm 1, \pm 5, \pm 9\})$, actually demonstrate that the one exception regarding types of admissible characters that we mentioned in Section~\ref{sec:eig0} can indeed occur. Namely, the admissible $5 \in \A_{60}(\{\pm 1, \pm 5, \pm 9\})$ is at the same time of type~1 and of type~2.

\section*{Acknowledgments}

\v S. Miklavi\v c acknowledges financial support by the Slovenian Research Agency (research program P1-0285 and research projects N1-0062, J1-1695, N1-0140, N1-0159, J1-2451, N1-0208, J3-3001, J3-3003).
P. \v Sparl acknowledges financial support by the Slovenian Research Agency (research program P1-0285 and research projects J1-1694, J1-1695, J1-2451, J1-3001).

\end{document}